\documentclass[12pt,a4paper]{amsart}
\usepackage[LGR,T1]{fontenc}
\usepackage[all]{xy}
\usepackage{graphicx}

\newenvironment{automates}
{
    \[
	\SelectTips{cm}{11}
    \UseTips
    \entrymodifiers={+[o][F]}
    \UseTips \xymatrix @+0.23mm }
{   \]}

\newtheorem{defi}{Definition }
\newtheorem{prop}[defi]{Proposition}
\newtheorem{rmk}[defi]{Remarque}

\newcommand{\Z}{\mathbb{Z}}

\begin{document}
\title{ On a property discovered by Xavier Grandsart}
\author{Pierre Arnoux}

\begin{abstract}We give a simple proof of the property discovered by Xavier Grandsart: let $W$ be a circular binary word; then the differences in the number of occurences $|W|_{0011}-|W|_{1100}$, $ |W|_{1101}-|W|_{1011}$ , $|W|_{1010}-|W|_{0101}$ and  $|W|_{0100}-|W|_{0010}$ are equal.
  \end{abstract}
\maketitle
  
 We are interested in the study of binary words, that is, finite sequences $W=W_0W_1\ldots W_{n-1}$ taking  values in the alphabet $\{0,1\}$; $n$ is the {\em length }
of the word $W$, denoted by $|W|$. We say that the word $U=U_0\ldots U_{k-1}$ of length $k\le n$ {\em occurs } in $W$ at position $i$ if $U_j=W_{i+j}$ for $0\le j\le k-1$, and in that case we say that $i$ is an {\em occurence } of $U$ in $W$, and that $U$ is a {\em factor} of $W$.
 
 To avoid special effects due to the extremities of the word $W$, we will consider {\em circular words} of length $n$, that is, words indexed by the cyclic group of integers mod $n$. This allow occurences at the end of the word; for example, the circular word  $W=0001$ admits the factor $010$, with an occurence at index 2.
 
 \begin{rmk} One can equivalently define a circular word as an infinite word which is periodic, of period W, and in that case we consider occurences as defined modulo $|W|$. If we identify two circular words which differ only in a shift of index, as 001 and 010, we can also define a circular word as a conjugacy class of words, where two words $W,W'$ are conjugate if there exist two words $U,V$ such that $W=UV$ and $W'=VU$. This is the usual definition in combinatorics of words. 
 \end{rmk}
 
We consider the number of occurences of a factor $U$ in a circular word $W$, which we will denote by $|W|_U$. For example, if $W=00101$, $|W|_{010}=2$, since $010$ occurs in position 1 and 3 in $W$. If we take a random circular word $W$ of large length and a fixed factor $U$, there is of course a probabilistic aspect 
to $|W|_U$, but there are also combinatorial relations. For example, it is easy to check that $|W|_{001}=|W|_{100}$; these two numbers are also the number of {\em runs} (sequences of maximal length) of 0 of length at least 2, since each such sequence starts with 100 and ends with 001. 

Around 2010, Xavier Grandsart, based on numerical experimentations, generalized this remark; he conjectured and proposed as a challenge the following proposition : 

\begin{prop} Let $W$ be a circular binary word. The differences in the number of occurences $|W|_{0011}-|W|_{1100}$, $ |W|_{1101}-|W|_{1011}$ , $|W|_{1010}-|W|_{0101}$ and  $|W|_{0100}-|W|_{0010}$ are equal.
\end{prop}

The challenge was quickly solved by 3 persons, Maher Younan\cite{Yo} who proposed a proof quite similar to the proof given here, Alberto Costa \cite {AC} and Pierre Deligne \cite{PD}, who proposed proofs by induction. We propose here an elementary proof based on De Bruijn graphs.

We start with some elementary remarks. Recall that a palindrome is a word $W$ which is equal to its mirror image, that is, if $n$ is the length of $W$, $W_i=W_{n-1-i}$ for all $i<n$. A palindromic pair is a pair $(U,V)$ of words of same length $n$ which are not palindromes and which are mirror image of each other.

There are 16 binary words of length 4; 4 of them, $1111$, $1001$, $0110$ and $0000$ are palindromes. 4 of them contain a run of length 3, and form two palindromic pairs whose two elements have same number of occurences, by the same proof as above: $(1000, 0001)$ and $(1110, 0111)$. The remaining 8 words form the 4 palindromic pairs considered by Xavier Grandsart: $(1010, 0101)$, $(0010, 0100)$, $(1011,1101)$, $(0011, 1100)$. These 8 words have the common property that their prefix of length 3 ends with two different letters: $001, 010, 101, 110$.

We will use the DeBruijn graph on binary words, and we recall the definition:

\begin{defi} The De Bruijn graph $B(2, n)$ is the oriented graph whose vertices are the binary words of length $n$, and whose edges are the binary words of length $n+1$; the initial vertex of the edge $U$ is the prefix of length $n$ of $U$, and the final vertex is the suffix of length $n$ of $U$.
\end{defi}

\begin{figure}[h]
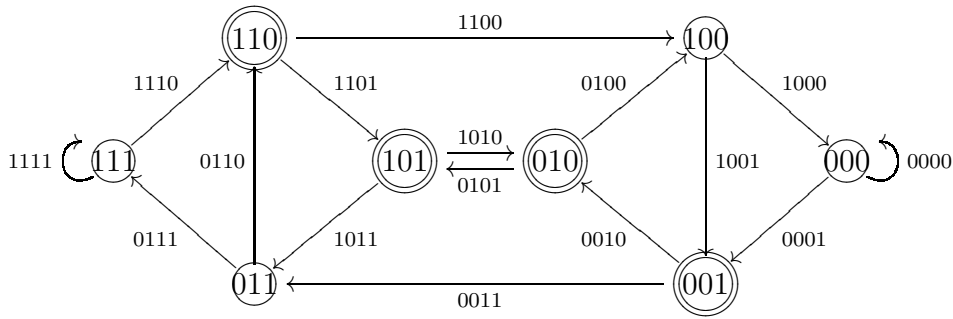

\begin{center}

\begin{automates}
{
*+\txt{}&*++[o][F=]{110}  \ar[rrr]^{1100}\ar[dr]^{1101}&   *+\txt{}&   *+\txt{}&100 \ar[dr]^{1000} \ar[dd]^{1001}  & *+\txt{}\\
\ar@(dl,ul)[]^{1111}111\ar[ur]^{1110}&*+\txt{}  &*++[o][F=] {101}\ar@<3pt>[r]^{1010} \ar[dl]^{1011}&*++[o][F=]{010} \ar@<3pt>[l]^{0101} \ar[ur]^{0100}&*+\txt{}  & 000  \ar@(dr,ur)[]_{0000} \ar[dl]^{0001}\\
*+\txt{}&011 \ar[ul]^{0111} \ar[uu]^{0110}&   *+\txt{}&   *+\txt{}&*++[o][F=]{001} \ar[lll]^{0011}  \ar[ul]^{0010} & *+\txt{}\\
}
\end{automates}
\caption{The De Bruijn graph $B(2,3)$}
\label{fig.B23}
\end{center}
\end{figure}

Figure \ref{fig.B23} shows the graph $B(2,3)$, were the initial vertices of the edges labeled by the words considered by Grandsart  are surrounded by double circles.
Each circular word of length $n$ defines a closed path in the De Bruijn graph $B(2,3)$; the sequence of consecutive vertices of this path is the sequence of factors of length 3 in their order of occurence, and the sequence of consecutive edges is the sequence of factors of length 4. 

\begin{proof}[First proof of the proposition]
Each vertex of the graph defines a relation in the number of occurences of factors of length 4. Indeed, if $W$ is a circular binary word, and $U$ is a factor of $W$ of length 3, then any occurence of $U$ is also an occurence of some $aU$ and some $Ub$, with $a,b\in \{0,1\}$. Hence we have $|W|_U=|W|_{U0}+|W|_{U1}=|W|_{0U}+|W|_{1U}$: the number of occurences of factor of length 4 satisfy Kirchhoff's law on the De Bruijn graph.

Taking $U=000$, we obtain $|W|_{0000}+|W|_{0001}=|W|_{0000}+|W|_{1000}$, hence $|W|_{0001}=|W|_{1000}$ as we saw above; similarly, using vertex $111$ we obtain $W|_{0111}=|W|_{1110}$. If we consider $U=101$ we get  $|W|_{101}= |W|_{1011}+|W|_{1010}=|W|_{0101}+|W|_{1101}$, and similarly considering vertex $010$, we have  $|W|_{1010}+|W|_{0010}=|W|_{0100}+|W|_{0101}$; these are equivalent forms for the relations $ |W|_{1101}-|W|_{1011}=|W|_{1010}-|W|_{0101}$ and  $|W|_{1010}-|W|_{0101}=|W|_{0100}-|W|_{0010}$.

The last relation is slightly more difficult to prove; the four remaining vertices give four relations, but they are not independent: the sum of the relations for all eight vertices is trivial, since each factor appears once on the right side and once on the left side of the equation. One can check that summing the relation for vertices $011$ and $110$, and using the fact that $|W|_{0111}=|W|_{1110}$, we recover the last relation in the proposition.
\end{proof}

We will give another proof which gives some meaning to the differences in the number occurences that appear in the proposition. 

\begin{proof}[Second proof of the proposition] If we start from the vertex $110$, we can either go directly to $101$ through the edge $1101$, or, starting with the edge $1100$, to $100$, then $n\ge 0$ times in $000$, then to $001$; a similar remark can be made for the other three vertices $101$, $010$, $001$, which are surrounded by double circles in the picture; each of them can only lead to two other vertices.

It follows that, if we erase, in the sequence of edges labelling the path of a circular word, the factors which are not in the 4 palindromic pairs considered by Grandsart, the remaining sequence labels a path in the smaller graph shown in Figure \ref{fig.short}.

\begin{figure}[!h]
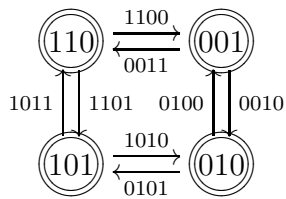

\begin{center}
\begin{automates}
{
*++[o][F=]{110}  \ar@<3pt>[r]^{1100}\ar@<3pt>[d]^{1101}&*++[o][F=]{001} \ar@<3pt>[d]^{0010} \ar@<3pt>[l]^{0011}  \\
*++[o][F=] {101}\ar@<3pt>[r]^{1010} \ar@<3pt>[u]^{1011}&*++[o][F=]{010} \ar@<3pt>[l]^{0101} \ar@<3pt>[u]^{0100}\\
}
\end{automates}
\caption{The shortened graph}
\label{fig.short}
\end{center}
\end{figure}

Between each pair of adjacent vertices on this graph, there are two opposite edges labelled by a palindromic pair. The graph is a cyclic graph; let us give this graph an orientation, so that edges labelled $0011, 1101, 1010$ and $0100$ correspond to a quarter turn in the positive direction, and the four other  edges go in the negative direction. 
A closed path on this graph makes a number $k$ of turns, counted with orientation, and for any pair of edges between adjacent vertices, the number of occurences of the positive edge minus the number of occurences of the negative edge is equal to $k$, hence it does not depend on the edge. This proves the proposition.

For a more formal version of this proof, define the permutation $R$ on the set of vertices of the square graph $001\to 110\to 101\to 010\to 001$; this is a cycle of order four. Any closed path on this graph is defined by a sequence $V_0, V_1, \ldots, V_n=V_0$ such that $V_{i+1}=R^{\epsilon_i}V_i$; it is also defined by the initial vertex $V_0$ and the sequence $R^{\epsilon_i}$, with $\epsilon_i=\pm 1$. We have $V_n= R^{\sum \epsilon_i} V_0= V_0$; since $R$ is a cycle of order 4, we must have $\sum_{i=0}^{n-1}\epsilon_i=4k$, for some $k\in \Z$. Any $R^{-1}R$ or $RR^{-1}$ in this coding corresponds to a pair of opposite edges labelled by a palindromic pair; simplifying this edge does not change the difference in the number of occurences of the palindromic pair, hence after simplification we can suppose that the path is coded by the initial vertex $V_0$, and  $R^{4k}$, so that each edge in the positive direction if $k>0$, or in the negative direction if $k<0$, occurs $|k|$ times, and the edges  in the other direction do not occur, hence the result.
\end{proof}

There is a simple way to compute the number $k$. Any letter $a$  occurs in a circular word either as an isolated letter $...bab...$ or as a run of length at least 2 $...aa...$. We can decompose a circular word in maximal sequences of isolated letters $...ababab...$ and sequences of runs of length at least 2; the number $k$ is equal to the difference between the number of maximal sequences of even length of isolated letters starting with 0, and the number of maximal sequences of even length of isolated letters starting with 1. In particular, for $W=010011$, we have $k=1$, and for $W=101100$, we have $k=-1$, as we can check on the graph.

The Kirchoff relations on the De Bruijn graph imply  that the 16  functions $W\mapsto |W|_U$, where $U$ is a binary word of length 4, satisfy 7 linear relations. Indeed, the 9 functions $W\mapsto |W|_U$, with $U=0000$ or $U=1V$ are independent (the corresponding words, if we remove the loops $0000$ and $1111$, generate a spanning tree for the graph), and the seven remaining functions, for $U$ starting with 0 and different from $0000$, can be generated from the other 9. 

Julien Cassaigne \cite{Ca} has pointed to me that this result can be widely generalized.  Consider an alphabet with $d$ letters, and  the functions $W\to |W|_U$, where $U$ is any word of length at most $l$ on this alphabet. There are $\sum_{k=0}^l d^k$ such functions; they are obviously not independent, since $|W|_U=|W|_{U0}+|W|_{U1}$, hence they are all generated by the $d^l$ functions defined by words of length $l$;  one proves that  the space generated by these functions has in fact dimension $(d-1)d^{l-1}+1$.  More precisely, it  can be deduced from \cite{CKS} that these functions can also be computed from the $W\mapsto |W|_U$, where $U$ are the words of length at most $l$ whose first and last letter is not 0. 

Remark that the function $W\mapsto |W|_U$, for $|U|<l$, is easily computed from the functions $W\mapsto |W|_V$,  for all words $V$ of length $l$, since it is the sum of all the occurences of the words of length $l$ which admit the word $U$ as prefix. These last functions must satisfy the Kirchhoff equations for the graph $B(d,l-1)$, so we recover the same number  $d^l-d^{l-1}+1$, cyclomatic number of the graph,  for the number of independent functions.

In our case, the sixteen functions generated by binary words of length 4 generate a space of dimension 9, and there are 7 independent relations, as we saw above.

\end{document}